# THE MINIMAL SEQUENCE OF REIDEMEISTER MOVES BRINGING THE DIAGRAM OF $(n+1, n)$-TORUS KNOT TO THAT OF $(n, n+1)$-TORUS KNOT

CHUICHIRO HAYASHI AND MIWA HAYASHI

ABSTRACT. Let $D(p, q)$ be the usual knot diagram of the $(p, q)$-torus knot, that is, $D(p, q)$ is the closure of the $p$-braid $(\sigma_1^{-1} \sigma_2^{-1} \cdots \sigma_{p-1}^{-1})^q$. As is well-known, $D(p, q)$ and $D(q, p)$ represent the same knot. It is shown that $D(n+1, n)$ can be deformed to $D(n, n+1)$ by a sequence of $\{(n-1)n(2n-1)/6\} + 1$ Reidemeister moves, which consists of a single RI move and $(n-1)n(2n-1)/6$ RIII moves. Using cowrithe, we show that this sequence is minimal over all sequences which bring $D(n+1, n)$ to $D(n, n+1)$.

## 1. INTRODUCTION

A Reidemeister move is a local move of a link diagram as in Figure 1. RI (resp. II) move creates or deletes a monogon face (resp. a bigon face). RIII move is performed on a 3-gon face, deleting it and creating a new one. Any such move does not change the link type. As Alexander and Briggs [1] and Reidemeister [6] showed, for any pair of diagrams $D_1$, $D_2$ which represent the same link type, there is a finite sequence of Reidemeister moves which deforms $D_1$ to $D_2$.

In [3], the knot diagram invariant cowrithe is introduced. It is not changed by an RI move, and is changed by zero or one by an RII move and by one by an RIII move. Hence, if two knot diagrams represent the same knot, then the difference of their cowrithes gives a lower bound for the number of RII and RIII moves required for deforming one to the other. In [2], Carter, Elhamdadi, Saito and Satoh gave a lower bound for the number of RIII moves by using extended n-colorings of knot diagrams in $\mathbb{R}^2$. Hass and Nowik introduced a certain knot diagram invariant by using smoothing and linking number in [4], and gave in [5] an example of an infinite sequence of diagrams of the trivial knot such that the $n$-th one has $7n - 1$ crossings, can be unknotted by $2n^2 + 3n$ Reidemeister moves, and needs at least $2n^2 + 3n - 2$ Reidemeister moves for being unknotted.

In this paper, we study the minimal number of Reidemeister moves on diagrams of torus knots, and show that the cowrithe gives exact lower bounds for $(n+1, n)$-torus knots. Let $D(p, q)$ be the usual knot diagram of the $(p, q)$-torus knot. See Figure 2, where $D(4, 3)$ is

The last author is partially supported by Grant-in-Aid for Scientific Research (No. 18540100), Ministry of Education, Science, Sports and Technology, Japan.

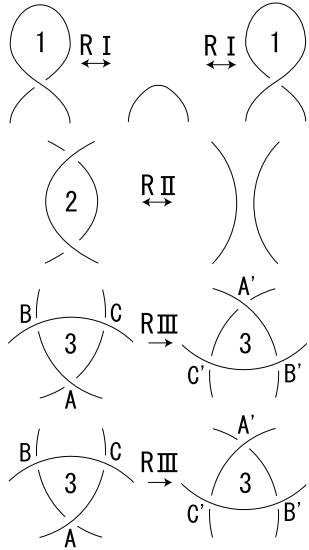

FIGURE 1

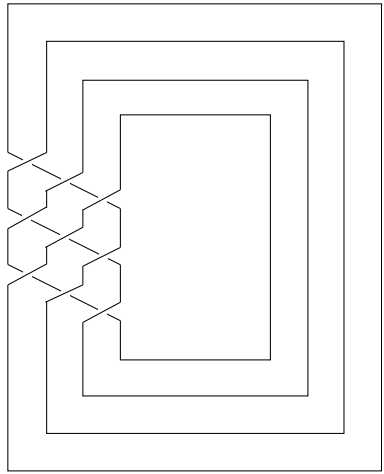

FIGURE 2

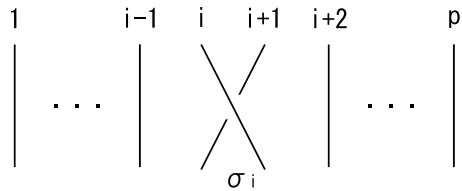

FIGURE 3

depicted. To be precise, for any $i \in \{1, 2, \cdots, p-1\}$, let $\sigma_i$ be the generator of the $p$-braid group $B_p$, which denotes the braid where the $i$-th strand crosses over the $(i+1)$ st strand (Figure 3). We let $b_i$ denote $\sigma_i^{-1}$ for short. Then, $D(p,q)$ is the closure of the $p$-braid $(b_1 b_2 \cdots b_{p-1})^q$. Thus $D(p,q)$ is a positive knot diagram. In fact, it has $(p-1)q$ positive crossings and no negative crossing with respect to the usual definition of sign of a crossing (Figure 4). As is well-known, $D(p,q)$ and $D(q,p)$ represent the same knot, and hence there is a sequence of Reidemsiter moves which brings $D(p,q)$ to $D(q,p)$. Using cowrithe, we obtain the following theorem.

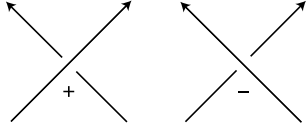

FIGURE 4

**Theorem 1.1.** *For any integer $n$ larger than or equal to 2, the usual positive knot diagram $D(n+1, n)$ of the $(n+1, n)$-torus knot can be deformed to $D(n, n+1)$ by a sequence of $(n-1)n(2n-1)/6+1$ Reidemeister moves, which consists of $(n-1)n(2n-1)/6$ RIII moves and a single RI move deleting a monogon face, and contains no move creating a negative crossing. Moreover, any sequence of Reidemeister moves bringing $D(n+1,n)$ to $D(n,n+1)$ must contains at least one RI move and at least $(n-1)n(2n-1)/6$ RII or RIII moves. Hence, the above sequence is minimal.*

The sequence of Reidemeister moves in the above theorem is described in Section 2. The estimation for the number of RI moves is easily obtained by using writhe. The definition of cowrithe is reviewed in Section 3, where we calculate cowrithes of $D(n+1,n)$ and $D(n,n+1)$. The proof of the theorem is given in Section 4. Notes for general $D(p,q)$ and positive knots are given in Section 5.

## 2. Deformation of $D(n+1, n)$ to $D(n, n+1)$

In this section, we deform $D(n+1,n)$ to $D(n,n+1)$ by a sequence of $(n-1)n(2n-1)/6+1$ Reidemeister moves. The deformation is described in Figures 5 through 12. First we move the overpath $\gamma_0$ going through the first $b_n$, the second $b_{n-1}$, $\cdots$ and the last $b_1$ of $D(n+1,n)$, so that $\gamma_0$ passes over all the crossings below it (Figure 5). This can be achieved by $(n-1)n/2$ RIII moves. Then we have a monogon face as in Figure 6, and it is eliminated by a single RI move. A closed $n$-braid is obtained as in Figure 7. When $n=2$, we have obtained $D(3,2)$.

When $n \geq 3$, we perform the deformations of $n-2$ steps as in Figures 7 through 11. The bold line is moved to the broken line in each step. The deformation of the $k$th step moves the $(n-k)$th strand of the braid except near the last crossing and is accomplished by $k \times (n-k-1) + (n-k-1) \times k = 2k(n-k-1)$ RIII moves. We begin the sequence of RIII moves with those on the trigonal faces marked with "3".

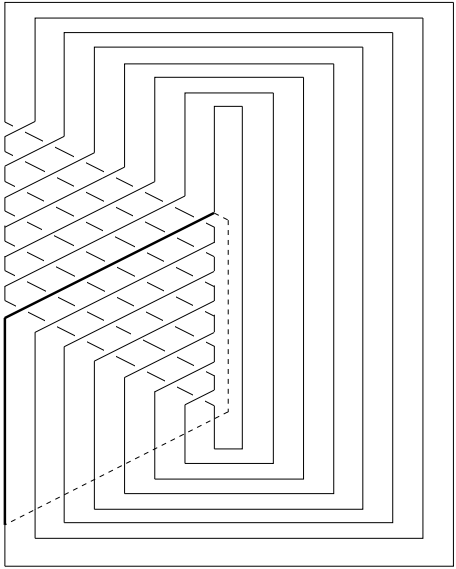

(8, 7)-torus knot

Figure 5

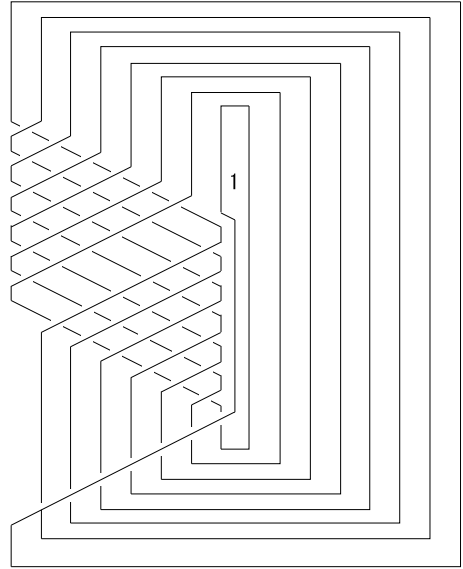

Figure 6

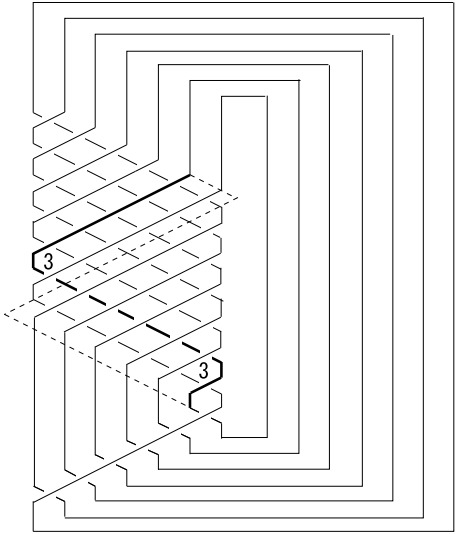

Figure 7

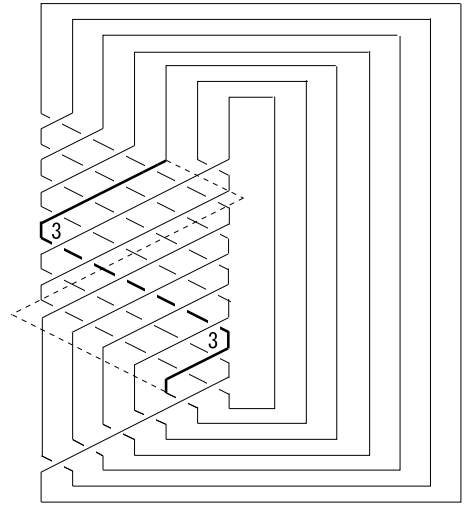

Figure 8

Thus we obtain the desired knot diagram $D(n, n+1)$ as in Figure 12, where the crossings within the triangle can be moved to the bottom of the braid by a Markov move.

The above deformation is composed of a single RI move, and $(n-1)n/2 + \sum_{k=1}^{n-2} 2k(n-k-1) = (n-1)n(2n-1)/6$ RIII moves. Thus the former half of Theorem 1.1 follows.

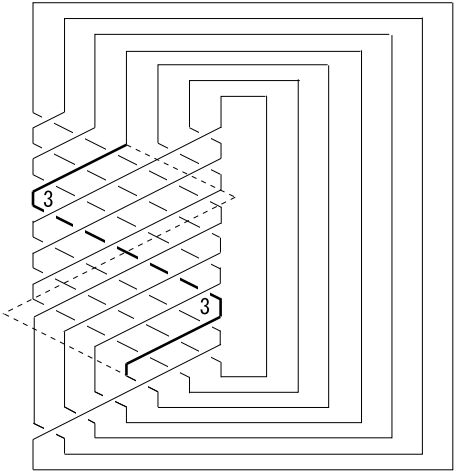

FIGURE 9

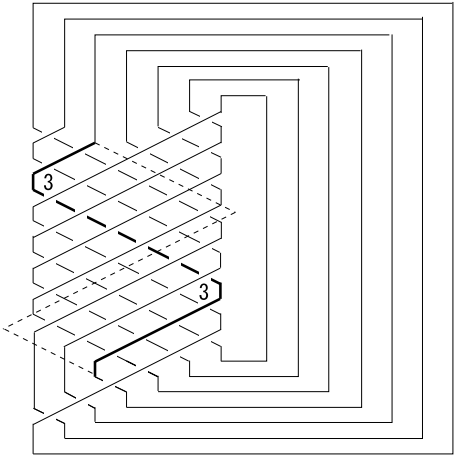

FIGURE 10

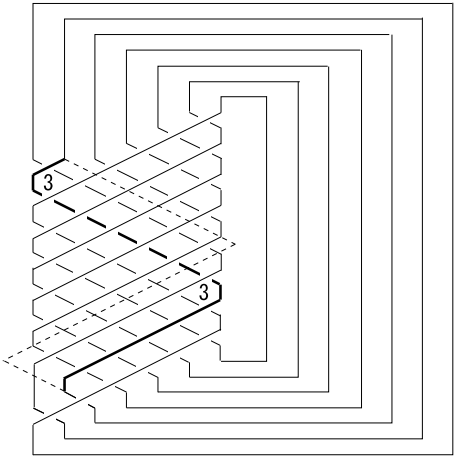

FIGURE 11

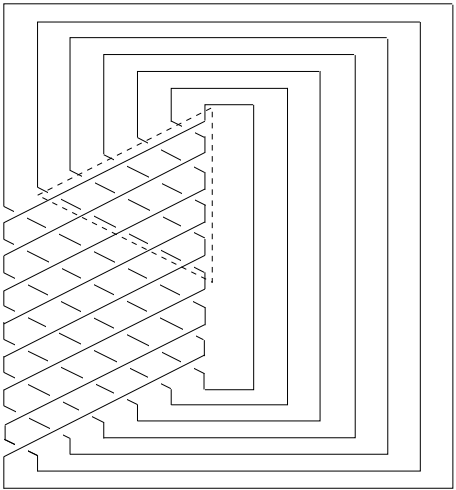

FIGURE 12

## 3. Calculation of cowrithe

In this section, we calculate cowrithes $x(D(n+1,n))$ and $x(D(n,n+1))$ of the diagrams $D(n+1,n)$ and $D(n,n+1)$.

**Lemma 3.1.** $x(D(n+1,n)) = (n-1)n^2(n+4)/6$ and $x(D(n,n+1)) = (n-1)n(n+1)^2/6$

Before beginning calculations, let us recall the definitions of chord diagram and cowrithe.

Let $D$ be a diagram in $S^2$ of a knot $K$. There are an embedding $f : S^1 \to S^3$ with $f(S^1) = K$, and a projection $\pi : S^3 - \{p_+, p_-\} \cong S^2 \times (-1,1) \to S^2 \times \{0\} \subset S^3$ with $\pi(K) = D$ where $p_+$ and $p_-$ are two points in $S^3$ disjoint from $K$. For every crossing $X$, let $X_1$ and $X_2$ be the preimage points of $X$, that is, $\pi(f(X_i)) = X$ for $i = 1$ and $2$. These two points $X_1$ and $X_2$ are contained in $S^1$, and $S^1$ is the unit circle in $\mathbb{R}^2$. The straight line segment connecting $X_1$ and $X_2$ in $\mathbb{R}^2$ is called the *chord* for the crossing $X$, and denoted by $\mathrm{Ch}_X$. The *chord diagram* $\mathrm{Ch}(D)$ of $D$ is the union of $S^1$ and the chords for all the crossings of $D$. For example, a knot diagram $6'_3$ and its chord diagram are described in Figure 13.

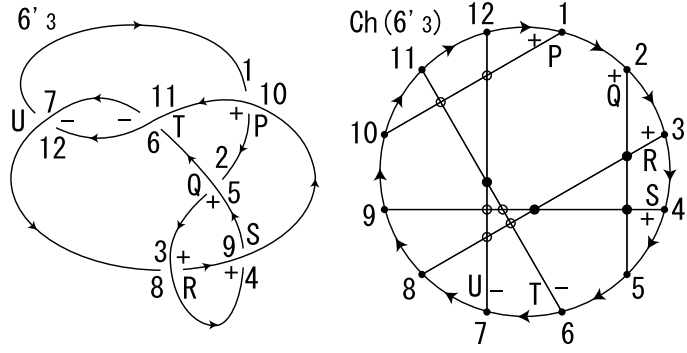

FIGURE 13

We say that two crossings of $D$ are *interleaved* if the chords for them intersects in a single point. A pair of crossings are interleaved if and only if their preimage points appear alternately on $S^1$. We give $D$ an arbitrary orientation. For every pair of interleaved crossings $P$ and $Q$ of $D$, we define the sign of the pair to be the product $(\mathrm{sign}\,P)\cdot(\mathrm{sign}\,Q)$. The *cowrithe* $x(D)$ of $D$ is the sum of signs of all the interleaved pairs of crossings of $D$. Reversing the orientation of $D$ does not change $x(D)$ since it does not change sign of any crossing. An example of calculation of cowrithe is described in the right of Figure 13, where we put a black dot at the intersection point of the chords $\mathrm{Ch}_U$ and $\mathrm{Ch}_T$ to indicate the positive sign of the pair of minus crosssings $U$ and $T$. The white dot at the intersection of the chords $\mathrm{Ch}_P$ and $\mathrm{Ch}_U$ stands for the negative sign of the pair $P$ and $U$. We can easily see that $x(6'_3) = 4 - 6 = -2$.

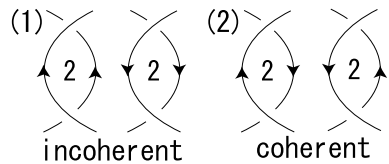

FIGURE 14

Cowrithe is not a knot invariant. However, it is almost an invariant in the sense that any Reidemeister move changes it at most by a constant.

**Theorem 3.2.** ([3]) *An RI move does not change $x(D)$. An RII move deleting a bigon face $f$ increases $x(D)$ by 1 if the orientations of the edges of $f$ are incoherent on the boundary circle of $f$ with respect to an orientation of $D$ as in Figure 14 (1). Otherwise, it does not change $x(D)$. An RIII move changes $x(D)$ by $\pm 1$.*

**Corollary 3.3.** *Let $D_1$ and $D_2$ be knot diagrams which represent the same knot. Then, any sequence of Reidemeister moves deforming $D_1$ to $D_2$ contains at least $|w(D_1) - w(D_2)|$ RII and RIII moves.*

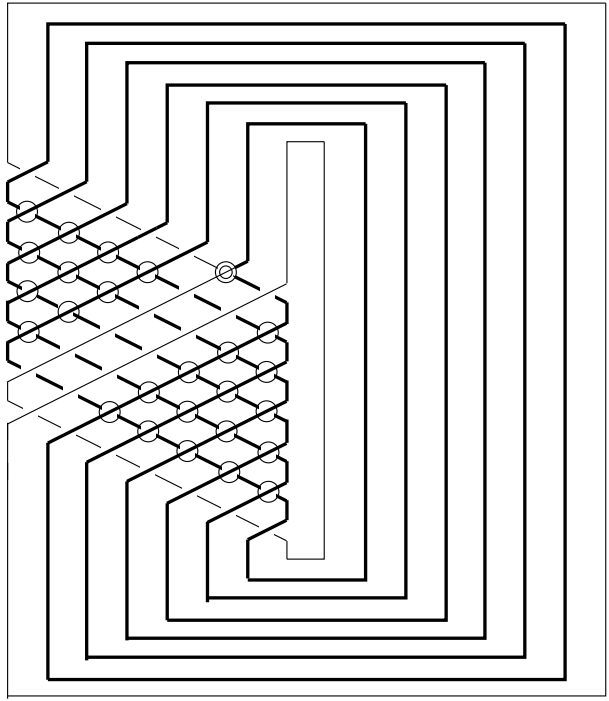

(8, 7)-torus knot

FIGURE 15

*Proof of Lemma 3.1.* First, we calculate $x(D(n+1, n))$.

Since the diagram $D(n+1, n)$ has only positive crossings, its cowrithe is equal to the number of pairs of interleaved crossings of $D$.

The preimage points of the first crossing point corresponding to $b_k$ ($1 \leq k \leq n$) divide the knot $K$ into two arcs. See Figure 15, where $b_k$ is the first $b_6$ doubly circled. One of the arcs, say $\lambda_k$, begins at the undercrossing point of the first $b_k$ and goes down along the first strand of the braid, then goes through the second, the third, $\cdots$, the $k$th strand and then goes back to the first $b_k$ at the overcrossing point. The arc $\lambda_k$ contains $2k(n-1)$ preimage points of crossings other than the endpoints. (Since $\lambda_k$ consists of $2k$ overpaths and underpaths, each of the two paths incident to the first $b_k$ contains $n-k$ preimage points and each of the other paths contains $n$ preimage points, the arc $\lambda_k$ contains $2(n-k) + 2(k-1)n = 2k(n-1)$ preimage points.) There are crossing points with both preimage points are on $\lambda_k$. $k-2$ $b_1$'s (the second through the $(k-1)$st), $k-3$ $b_2$'s (the second through the $(k-2)$nd), $\cdots$ and one $b_{k-2}$ (the second), one $b_{n-k+2}$ (the $k$th), two $b_{n-k+3}$'s (the $(k-1)$st and the $k$th), $\cdots$ and $(k-1)$ $b_n$'s (the second through $k$ th). Each of these $(k-2)(k-1)/2+(k-1)k/2 = (k-1)^2$ crossings and the first $b_k$ are not interleaved. See Figure 15, where such crossings are circled. Hence the first $b_k$ contributes to $x(D(n+1, n))$ by $2k(n-1) - 2(k-1)^2$. Any crossing $b_k$ has the same contribution by symmetry of the diagram $D(n+1, n)$.

The diagram $D(n+1, n)$ has $n$ crossings corresponding to $b_k$ with $1 \le k \le n$. We counted contributions doubly. Hence
$x(D(n+1, n)) = (1/2) \sum_{k=1}^{n} \{(2k(n-1) - 2(k-1)^2) \times n\}$
$= \{(n-1)n(n+1)/2 - (1/6)(n-1)n(2(n-1)+1)\} \times n = (n-1)n^2(n+4)/6.$

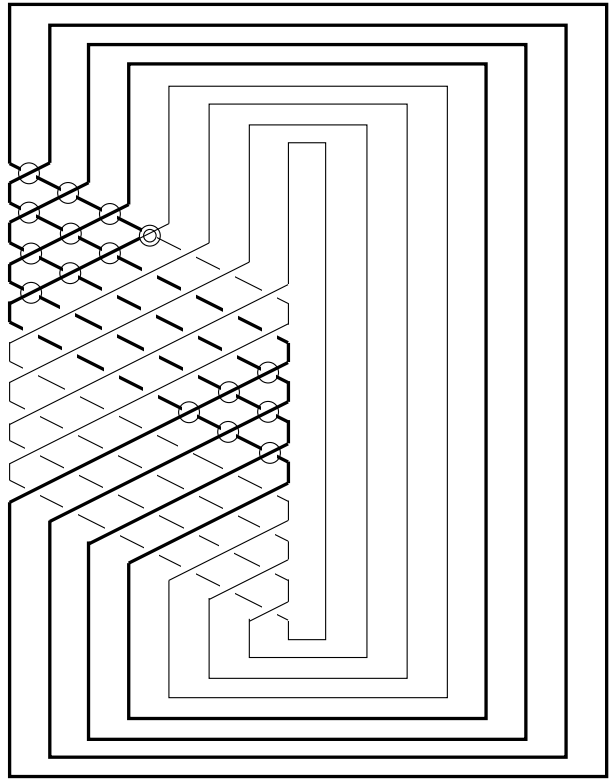

(8, 9)-torus knot

FIGURE 16

Next, we calculate $x(D(n, n+1))$. The preimage points of the first crossing point corresponding to $b_k$ $(1 \le k \le n-1)$ divide the knot $K$ into two arcs. See Figure 16, where $b_k$ is the first $b_4$ doubly circled. One of the arcs, say $\mu_k$, begins at the undercrossing point of the first $b_k$ and goes up along the first strand of the braid, then goes through the second, the third, $\cdots$, the $(k+1)$ st strand, to go back to the first $b_k$ at the overcrossing point. The arc $\mu_k$ contains $2(kn-1)$ preimage points of crossings other than the endpoints. (Since $\mu_k$ consists of $2(k+1)$ overpaths and underpaths, each of the two paths incident to the first $b_k$ contains $k-1$ preimage points and each of the other paths contains $n-1$ preimage points, the arc $\mu_k$ contains $2(k-1) + 2k(n-1) = 2(kn-1)$ preimage points.) There are crossing points with both preimage points are on $\mu_k$. $k$ $b_1$'s (the first through the $k$th), $k-1$ $b_2$'s (the first through the $(k-1)$st), $\cdots$ and two $b_{k-1}$ (the first and the second), one $b_{n-k+1}$ (the $(k+1)$st), two $b_{n-k+2}$'s (the $k$th and the $(k+1)$st), $\cdots$ and $(k-1)$ $b_{n-1}$'s (the third through $(k+1)$st). Each of these $k(k+1)/2 - 1 + (k-1)k/2 = k^2 - 1$ crossings and the first $b_k$ are not interleaved. See Figure 16, where such crossings are circled. Hence the first $b_k$ contributes to $x(D(n, n+1))$ by $2(kn-1) - 2(k^2-1)$. Any crossing $b_k$ has the same contribution by symmetry of the diagram $D(n, n+1)$.

The diagram $D(n, n+1)$ has $n+1$ crossings corresponding to $b_k$ with $1 \le k \le n-1$. Hence, we have

$x(D(n, n+1)) = (1/2) \sum_{k=1}^{n-1} \{(2(kn-1) - 2(k^2-1)) \times (n+1)\} = (n-1)n(n+1)^2/6$. □

## 4. Proof of Theorem 1.1

The former half of the theorem has been shown in Section 2.

We can obtain easily the estimation for the number of RI moves by using writhe. For any knot diagram $D$, the writhe $w(D)$ is the sum of the sign of all the crossings of $D$ with respect to an arbitrary orientation of $D$. Writhe is changed by an RI move by one, and is unchanged by an RII or RIII move. Note that $w(D)$ is unchanged if the orientation of $D$ is reversed. Hence, for two knot diagrams $D_1$ and $D_2$, we need at least $|w(D_1) - w(D_2)|$ RI moves to deform $D_1$ to $D_2$. Because $D(n+1, n)$ has $n^2$ positive crossings and no negative crossings, we have $w(D(n+1, n)) = n^2$. Similarly, $w(D(n, n+1)) = (n-1)(n+1) = n^2 - 1$, and hence $|w(D(n+1, n)) - w(D(n, n+1))| = 1$. Thus we need at least one RI move to deform $D(n+1, n)$ to $D(n, n+1)$.

By Lemma 3.1, the difference of cowrithes is $x(D(n+1, n)) - x(D(n, n+1)) = (n-1)n^2(n+4)/6 - (n-1)n(n+1)^2/6 = (n-1)n(2n-1)/6$. Hence, by Corollary 3.3, we need at least $(n-1)n(2n-1)/6$ RII and RIII moves to deform $D(n+1, n)$ to $D(n, n+1)$. This completes the proof of Theorem 1.1.

## 5. Some notes for generalization

It is very easy to see that $x(D(5, 2)) = 12$ and $x(D(2, 5)) = 10$. However, it seems to be impossible to deform $D(5, 2)$ to $D(2, 5)$ by a sequence of Reidemeister moves containing at most two RII and RIII moves.

It is very beautiful that two positive link diagrams representing the same link are connected by a sequence of Reidemeister moves which do not create a negative crossing. Note that such a sequence does not contain an RII move. However, the flype operation as shown in Figure 17 does not seem to be realizable by a sequence of Reidemeister moves without an RII move.

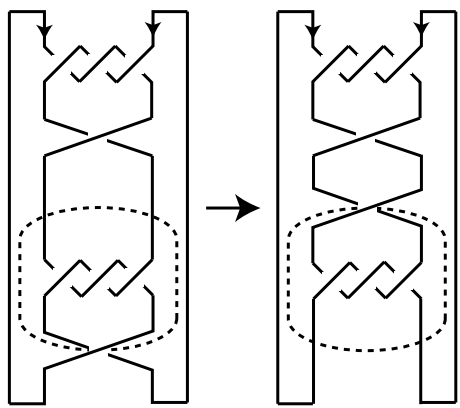

Figure 17

Chuichiro Hayashi: Department of Mathematical and Physical Sciences, Faculty of Science, Japan Women's University, 2-8-1 Mejirodai, Bunkyo-ku, Tokyo, 112-8681, Japan. hayashic@fc.jwu.ac.jp

Miwa Hayashi: Department of Mathematical and Physical Sciences, Faculty of Science, Japan Women's University, 2-8-1 Mejirodai, Bunkyo-ku, Tokyo, 112-8681, Japan. miwakura@fc.jwu.ac.jp